\documentclass[a4paper,11pt]{article}

\pdfoutput=1
\usepackage{amssymb,amsfonts,framed,xcolor,graphicx}
\usepackage[colorlinks=true]{hyperref}

\newtheorem{thm}{Theorem}[section]
\newtheorem{dfn}[thm]{Definition}

\newtheorem{cor}[thm]{Corollary}
\newtheorem{con}[thm]{Conjecture}
\definecolor{shadecolor}{gray}{0.95}

\begin{document}

\begin{center}
{\bf\large THE COVERING RADIUS PROBLEM FOR SETS OF PERFECT MATCHINGS}\bigskip

Cheng Yeaw \underline{Ku}\medskip

{\it Department of Mathematics, National University of Singapore, S117543}\bigskip\bigskip

Alan J. \underline{Aw}\medskip

{\it Raffles Institution (Junior College), One Raffles Institution Lane, S575954}\bigskip\bigskip

\begin{abstract}

Consider the family of all perfect matchings of the complete graph $K_{2n}$ with $2n$ vertices. Given any collection $\mathcal M$ of perfect matchings of size $s$, there exists a maximum number $f(n,x)$ such that if $s\leq f(n,x)$, then there exists a perfect matching that agrees with each perfect matching in $\mathcal M$ in at most $x-1$ edges. We use probabilistic arguments to give several lower bounds for $f(n,x)$. We also apply the Lov\'{a}sz local lemma to find a function $g(n,x)$ such that if each edge appears at most $g(n, x)$ times then there exists a perfect matching that agrees with each perfect matching in $\mathcal M$ in at most $x-1$ edges. This is an analogue of an extremal result vis-\`{a}-vis the covering radius of sets of permutations, which was studied by Cameron and Wanless (cf. \cite{cameron}), and Keevash and Ku (cf. \cite{ku}). We also conclude with a conjecture of a more general problem in hypergraph matchings. 
\end{abstract}

\end{center}

\bigskip



\section{Introduction}

In this paper, let $K_{2n}$ be the complete graph with $2n$ {\it vertices}, $n\in\mathbb N$. A {\it matching} in $K_{2n}$ is a set of pairwise non-adjacent {\it edges}; that is, no two edges share a common vertex. A {\it perfect matching} is a matching which matches all vertices of the graph; that is, every vertex of the graph is incident to exactly one edge of the matching. Any perfect matching is represented by a collection of two-element sets where the elements of each set are two distinct vertices; for instance $\{\{v_1,v_2\},\{v_3,v_4\}\}$ is a perfect matching in $K_4$, as shown in figure 1. Also, a {\it hypergraph} is a pair $(V,E)$, where $V$ is a finite set of vertices and $E$ is a finite family of subsets of $V$, called {\it hyperedges}. Using the terminology in \cite{bollobas}, we denote by $H$ a $t${\it-uniform} hypergraph ($t\in\mathbb N$), a hypergraph with $E\subset V^{(t)}$, where $V^{(t)}:=\{Y:Y\subset V,|Y|=t\}$.\bigskip

\begin{figure}[h]
	\centering
		\includegraphics[width=0.50\textwidth]{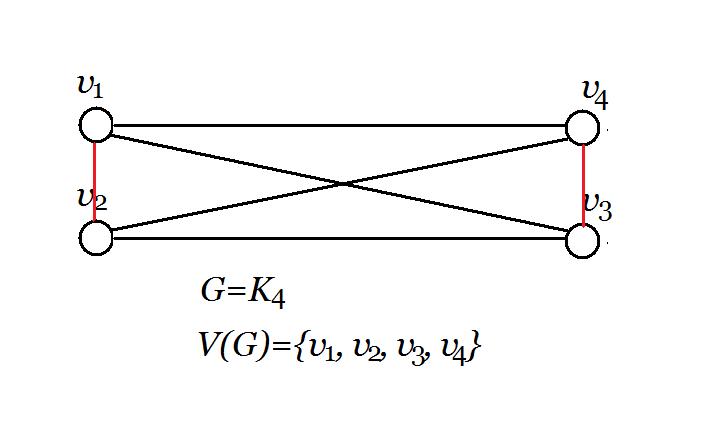}
	\caption{Perfect matching in $K_4$.}
	\label{fig:k4}
\end{figure}

Consider the following problem:\bigskip 

\begin{framed}
Given a collection of similar structures (family of permutations, system of finite sets et al.), what is the maximum size of the collection in order to ensure the existence of another such structure that shares at most $k$ elements with each structure in the collection? 
\end{framed}

\bigskip

This problem is also known as the covering radius problem. This research problem has its origins in group theory, particularly in the permutation group $S_n$ acting on the set $[n]$ of natural numbers from $1$ to $n$. In any collection $G$ of permutations, we can measure the {\it Hamming distance} (or distance) $d(g,h)$ between a permutation $g$ in $G$ and any permutation $h$ picked from $S_n$. Here, the Hamming distance between two permutations is the number of positions in which they defer. For example, in $S_3$, $d(123,231)=3$. If we were to fix $h$ above and measure the distances $d(h,p)$ for every $p\in G$, there exists a minimum distance which we can obtain between $h$ and some\footnote{There may exist more than one choice of $p_0$ which gives a minimum distance.} $p_0\in G$, i.e. $\min\{d(h,p):p\in G\}=d(g,p_0)$. Now, repeating this procedure for every permutation $h\in S_n$, we can find the maximum of all the minimum distances measured earlier. This maximum value, denoted $cr(G)$, is the covering radius of the collection $G$; in fact, a simple argument shows that $cr(G):=\max_{h\in S_n}\min_{g\in G}{d(g,h)}$. Therefore, the covering radius problem is the problem of finding or, in many cases, estimating the covering radius of any given collection of permutations.\bigskip

Apropos of recent research, lower bounds of the covering radius of $G$ of $S_n$ have been established by Cameron and Wanless (2005) in \cite{cameron}, in which covering arguments were used to formulate a general criteria to find lower bounds of covering radii of sets of permutations. Keevash and Ku (2006) later improved (cf. \cite{ku}) the general criteria, obtaining an even stronger result to determine the lower bound of the covering radius of any collection of permutations with some constraints vis-\`{a}-vis a frequency parameter. Their result is the best possible so far. Covering radius results have profound applications and implications in group theory and combinatorial structures; for instance, the authors above have applied their results to Latin squares and Latin transversals, and recent literature suggests several generalizations of this theory to general groups \cite{ku}. Moreover, Similar classes of problems for intersecting families of finite sets have been studied extensively, and in particular \cite{babai} and \cite{graham} are good sources of information.\bigskip

In this paper, we consider the analogue of the problem mentioned above for perfect matchings in complete graphs. Our fundamental question is as follows:\bigskip

\begin{framed}
Suppose we have a finite collection of perfect matchings of $K_{2n}$. What is the largest possible number of elements in this collection such that we can find a perfect matching of $K_{2n}$ that {\it agrees} with each perfect matching in the collection in at most $x-1$ edges?
\end{framed}

\bigskip

Let $M$ denote an arbitrary perfect matching of $K_{2n}$. Moreover, call a vertex set $U\subset V(K_{2n})$ {\it{good}} (w.r.t. $M$) if $|U|\leq n$ and $\forall v_i,v_j\in U(e_{ij}\notin M)$. For a set $W$ of vertices, we say that two perfect matchings $M,M'$ agree on $W$ iff (i) $W$ is good w.r.t. $M$ and $M'$; (ii) $\forall v_i\in W(e_{ij}\in M\Leftrightarrow e_{ij}\in M')$. For example, if $M=\{\{1,2\}, \{3,4\}, \{5,6\}\}$ and $M'=\{\{1,2\}, \{3,5\}, \{4,6\}\}$, then $M$ and $M'$ agree on $W=\{1\}$. We now present a few elementary bounds which are obtained from Boole's inequality (also known as union bound).


\section{Elementary Results}

The union bound gives us several useful results.

\begin{thm}
\label{Crude 1}
Let $\mathcal M=\{M_1,...,M_s\}$ be a collection of perfect matchings in $K_{2n}$. If $s\leq{x!}\cdot\frac{{2n-x\choose x}}{{n\choose x}}$, then there exists a perfect matching that agrees with each $M_i\in\mathcal M$ in at most $x-1$ edges.
\end{thm}

\begin{shaded}
\noindent\textbf{Proof.}\bigskip

Randomly select a perfect matching out of all the perfect matchings. Consider any $M_i\in\mathcal M$ and any $T\subset V$, $|T|=x$, which is good w.r.t. $M_i$. Let $A_{i,T}$ be the event that the perfect matching selected agrees with $M_i$ on $T$. Then $$\mathbb P(A_{i,T})=\left[\frac{(2n-2x)!}{2^{n-x}\cdot(n-x)!}\right]\bigg/\left[\frac{(2n)!}{2^n\cdot n!}\right],$$\bigskip

since there are $\frac{(2n)!}{2^n\cdot n!}$ perfect matchings and exactly $\frac{(2n-2x)!}{2^{n-x}\cdot(n-x)!}$ of them with $x$ fixed edges.\bigskip

Let us sum the probabilities over all possible $i$ and $T$. Clearly, there are $s$ possible values of $i$, and for each perfect matching $M_i$ there are less than ${2n\choose x}$ good $T$. Moreover, $A_{i,T}\Leftrightarrow A_{i,T'}$ whenever $T,T'$ belong the same collection of edges, implying that the number of good $T$ should be reduced by a factor of $2^x$ in our calculation in order to avoid counting same events more than once. Thus we have, by Boole's inequality, $$\mathbb P\left(\bigcup_{i,T}{A_{i,T}}\right)\leq\sum_{i,T}{\mathbb P(A_{i,T})}<s\cdot{2n\choose x}2^{-x}\cdot\left[\frac{(2n-2x)!}{2^{n-x}\cdot(n-x)!}\right]\bigg/\left[\frac{(2n)!}{2^n\cdot n!}\right]\leq 1.$$\bigskip

Therefore, with positive probability none of the $A_{i,T}$ occur, and there must exist a perfect matching which agrees with each $M_i\in\mathcal M$ on a vertex set of at most $x-1$ vertices, i.e. in at most $x-1$ edges.\smallskip

\begin{flushright}
$\blacksquare$
\end{flushright}
\end{shaded}

\bigskip

Notice that the bound on $s$ is weak, since we gave a crude bound of ${2n\choose x}2^{-x}$ for the number of good $T$. By considering $T$ differently, we can obtain the exact number of good $T$. Here, we introduce the notion of $x$-matchings.\bigskip

\begin{dfn}
\label{matching}
An $x$-matching of $K_{2n}$ is a matching of size $x$. Thus, if $x=n$, then the $x$-matching is simply a perfect matching.
\end{dfn}

With this in mind, we can derive a larger upper bound on $s$ as follows:\bigskip

\begin{thm}
\label{Crude 2}
Let $\mathcal M=\{M_1,...,M_s\}$ be a family of perfect matchings in $K_{2n}$. If $s<\frac{x!}{2^x}\cdot\frac{{2n\choose x}\cdot{2n-x\choose x}}{{n\choose x}^2}$, then there exists a perfect matching that agrees with each $M_i\in\mathcal M$ in at most $x-1$ edges.
\end{thm}

\bigskip

\begin{shaded}
\noindent\textbf{Proof.}\bigskip

Randomly pick a perfect matching out of all the perfect matchings. Consider any $M_i\in\mathcal M$ and pick any $x$-matching $X\subset M_i$. Let $A_{i,X}$ be the event that the perfect matching picked contains $X$. Then $$\mathbb P(A_{i,X})=\left[\frac{(2n-2x)!}{2^{n-x}\cdot(n-x)!}\right]\bigg/\left[\frac{(2n)!}{2^n\cdot n!}\right],$$\bigskip

by the same reasoning as shown in the proof of Theorem \ref{Crude 1}.\bigskip

Let us sum the probabilities over all $i$ and $X$. Clearly, there are $s$ possible values of $i$, and for each perfect matching $M_i$ there are exactly ${n\choose x}$ $x$-matchings. Thus we have $$\mathbb P\left(\bigcup_{i,X}{A_{i,X}}\right)\leq\sum_{i,X}{\mathbb P(A_{i,X})}\leq s\cdot{n\choose x}\cdot\left[\frac{(2n-2x)!}{2^{n-x}\cdot(n-x)!}\right]\bigg/\left[\frac{(2n)!}{2^n\cdot n!}\right]< 1.$$\bigskip

Therefore, with positive probability none of the $A_{i,X}$ occur, and there must exist a perfect matching which agrees with each $M_i\in\mathcal M$ in at most $x-1$ edges.\smallskip

\begin{flushright}
$\blacksquare$
\end{flushright}
\end{shaded}


\section{Main Result}

Our main theorem requires the Lov\'{a}sz sieve. The Lov\'{a}sz local lemma is a powerful tool for showing the existence of structures with desired properties. Briefly speaking, we toss our events onto a probability space and evaluate the conditional probabilities of certain bad events occurring. If these probabilities are not too large in value, then with positive probability none of the bad events occur. More precisely,\bigskip

\begin{thm}[Lov\'{a}sz]
\label{Lovasz}
Let $\mathcal A=\{A_1,A_2,...,A_n\}$ be a collection of events in an arbitrary probability space $(\Omega,\mathcal F,\mathbb P)$. A graph $G(V,E)$ is called a dependency graph $(V=\{1,2,...,n\})$ for the events $A_1,A_2,...,A_n$, where $e_{ij}\in E$ iff $A_i$ and $A_j$ are related by some property $\pi$. Suppose that $G(V,E)$ is a dependency graph for the above events and $\exists x_1,x_2,...x_n\in[0,1)$, $S\subset\{1,2,...,n\}\setminus\{j:e_{ij}\in E\}$ such that $$\mathbb P\left(A_i\mid\bigcap_{k\in S}{\overline{A_k}}\right)\leq x_i\cdot\prod_{e_{ij}\in E}{(1-x_j)}.$$\bigskip

\noindent Then $\mathbb P(\bigcap_{i=1}^{n}{\overline{A_i}})\geq\prod_{i=1}^{n}{(1-x_i)}$. Equivalently, with positive probability none of the $A_i$ occur.
\end{thm}

\bigskip

The proof of Theorem \ref{Lovasz} can be found in chapter 5 of \cite{alon} and chapter 19 of \cite{jukna}. We used the following special case (cf. \cite{alon}) of Theorem \ref{Lovasz} in our result:\bigskip

\begin{cor}
\label{Symmetric}
Suppose that $\mathcal A=\{A_1,A_2,...,A_n\}$ is a collection of events, and for any $A_i\in\mathcal A$ there is a subset $\mathcal D_{A_i}\subset\mathcal A$ of size at most $d$, such that for any subset $\mathcal S\subset\mathcal A\setminus\mathcal D_{A_i}$ we have $\mathbb P\left(A_i\mid\bigcap_{A_j\in\mathcal S}{\overline{A_j}}\right)\leq p$. If $ep(d+1)\leq 1$, then $\mathbb P(\bigcap_{i=1}^{n}{\overline{A_i}})>0$.
\end{cor}

\bigskip

We now establish our main result on the covering radius problem for sets of perfect matching using Corollary \ref{Symmetric}. In this proof, the strategy we use mirrors that of the proof of the lower bound on the covering radius for sets of permutations, as presented in \cite{ku}. Such a strategy is also used in proof of the Erd\H{o}s-Spencer theorem on Latin transversals, as presented in chapter 5 of \cite{alon} (pp. 73-74).\bigskip

\begin{thm}
\label{Main}
Let $\mathcal M=\{M_1,...,M_s\}$ be a collection of perfect matchings in $K_{2n}$. Moreover, each of the ${2n\choose 2}$ edges appears at most $k$ times in the $M_i\in\mathcal M$ (we call $k$ the frequency parameter). If $k\leq\frac{1}{e\cdot2x(2n-1){n-1\choose x-1}}\left(\sum_{j=2x-n}^{x}{{2x\choose 2j}\frac{(2j)!}{j!\cdot2^j}}-e\right)$, then there exists a perfect matching which agrees with each perfect matching $M_i\in\mathcal M$ in at most $x-1$ edges.
\end{thm}

\begin{shaded}
\noindent\textbf{Proof.}\bigskip

Randomly pick a perfect matching $M$ from the set of all perfect matchings in $K_{2n}$. Consider any $M_i\in\mathcal M$ and any $x$-matching $X\subset M_i$. Let $A_{i,X}$ be the event that $X\subset M$. In our dependency graph, connect $A_{i,X}$ to $A_{i',X'}$ iff $X$ and $X'$ share at least one common vertex in their underlying vertex sets.\bigskip

For each $A_{i,X}$, let the set of its neighbours in the dependency graph be $\mathcal D_{i,X}$.\bigskip

\noindent{\it Claim: $|\mathcal D_{i,X}|\leq k\cdot2x(2n-1){n-1\choose x-1}=d$.}\bigskip

Indeed, we first pick one vertex out of the $2x$ vertices in $X$, then choose out of the $2n-1$ remaining vertices one particular vertex to be its neighbour in the perfect matching. Next, we choose a perfect matching $M_{i'}\in\mathcal M$ that contains the constructed edge; this can be done in at most $k$ ways. Lastly, we just pick $x-1$ out of the remaining $n-1$ edges in $M_{i'}$ to form $X'$.\bigskip

Let us now consider the probability $\mathbb P\left(A_{i,X}\mid\bigcap_{A_{i',X'}\in\mathcal S}{\overline{A_{i',X'}}}\right)=p_0$ for any subset $\mathcal S\subset \mathcal A\setminus\mathcal D_{i,X}$. For brevity let us label $E=\bigcap_{A_{i',X'}\in\mathcal S}{\overline{A_{i',X'}}}$. We shall bound $p_0$ from above.\bigskip

Fix $A_{i,X}$. Without loss of generality, let the underlying set of $2x$ vertices of the $x$-matching $X$ be $V=\{v_1,v_2,...,v_{2x}\}$. Now, randomly pair arbitrarily many of the $2x$ vertices in $V$. This gives us a collection of singletons (vertices) and doubletons (edges) - we call such a collection $W$. However, we restrict our $W$ such that the total number of singletons and doubletons in any $W$ cannot exceed $n$, i.e. if there are $2p$ singletons (the number of singletons must be even) and $q$ doubletons then $2p+q\leq n$ (note that $q+p=x$). This gives us $n-x \geq p \geq0$. Thus, the set $\mathcal W$ of all such restricted $W$ has cardinality $$\sum_{k=2x-n}^{x}{{2x\choose 2k}\frac{(2k)!}{k!\cdot2^k}},$$\bigskip

where each summand is the number of ways to partition the underlying vertex set into a collection of $k$ doubletons and $2x-2k$ singletons. (To resolve ambiguity in the expression above, let ${n\choose k}=0$ if $k<0$.)\bigskip

For each $W$, let $B_W$ be the event that a perfect matching contains $W$. For example, $$M=\{\{v_1,v_2\},\{v_3,v_4\},\{v_5,v_6\}\}$$ 

contains $W=\{\{v_1,v_2\},\{v_3\},\{v_6\}\}$, where every pair of singletons in $W$ does not belong to any edge in $M$. Clearly, all the $B_W$ are mutually exclusive, and their union equals $\Omega$.\bigskip

We shall show that the number of perfect matchings contained in $B_W\cap E$ is at least the number of perfect matchings contained in $A_{i,X}\cap E$. This can be done by means of constructing an injection from $A_{i,X}\cap E$ to $B_W\cap E$ for a particular fixed $W$.\bigskip

\noindent{\it Claim: $\forall W\in\mathcal W\left(|B_W\cap E|\geq |A_{i,X}\cap E|\right)$.}\bigskip

First, for any $M\in A_{i,X}\cap E$, consider the remaining $n-x$ edges not in $X$. Direct each edge such that the tail of the directed edge is the vertex with the smaller subscript. Thus, every edge $\{v_i,v_j\}\notin X, i<j$ becomes $(v_i,v_j)$. This gives us $n-x$ ordered pairs of vertices. Now, arrange the $n-x$ edges lexicographically by the following rule: compare every two edges and place the edge whose first component has a vertex with a smaller subscript in front; i.e. if $(v_i,v_j),(v_k,v_l)$ are both directed edges originally belonging to $W$ and $k<i$, then $(v_k,v_l)$ goes in front of $(v_i,v_j)$. This gives an ordered $(n-x)$-tuple $((v_a,v_b),(v_c,v_d),...)$ where $a<b$, $c<d$ and $a<c$ and so on. Denote this sequence of transformations on $M\setminus X$ by $\tau$. Notice that for any two distinct perfect matchings $M,M'\in A_{i,X}\cap E$, at least two of their $n-x$ edges outside $X$ are distinct (e.g. $\{\{v_\alpha,v_\beta\},\{v_\chi,v_\delta\}\}\subset M$, $\{\{v_\alpha,v_\chi\},\{v_\beta,v_{\delta'}\}\}\subset M'$), so their images under $\tau$ will also be distinct. Therefore, $\tau$ is injective.\bigskip

Now consider $W$. Without loss of generality, let $W$ contain $2p$ singletons, where $0\leq p\leq n-x$. Order the singletons in $W$ naturally by comparing their respective vertices' subscripts. Denote by $W_\gamma$ the image of $W$ under the natural ordering. Without loss of generality, write $$W_\gamma=\Big\{\big(\{v_1\},\{v_2\},...,\{v_{2p}\}\big),\{v_{2p+1},v_{2p+2}\},...,\{v_{2x-1},v_{2x}\}\Big\}.$$

We define a mapping as follows:\bigskip

For any $M\in A_{i,X}\cap E$, we consider their images under $\tau$. Treating the $(n-x)$-tuple, of which each component is an ordered pair, as an ordered string of vertices of length $2n-2x$, we select the first $2p$ vertices appearing in the string and pair the $k$th vertex in the string with $v_k$ in $W$. This gives us the set $\Gamma=\{(\{v_1,v_a\},\{v_2,v_b\},...),\{v_{2p+1},v_{2p+2}\},...\}$. Remove the natural ordering on $\Gamma$ to yield $\Gamma_0=\{\{v_1,v_a\},\{v_2,v_b\},...,\{v_{2p+1},v_{2p+2}\},...\}$. Now map the shortened string of length $2n-2x-2p$ back to its set of unordered $n-x-p$ edges (note that this gives us edges which were originally in $M$); call this edge set $\Gamma_1$. Clearly, $\mathring{\Gamma}(M)=\Gamma_0\cup\Gamma_1$ gives us a perfect matching in $B_W\cap E$, since $E$ is the event that $X'\not \subseteq M$ where $X'\cap X=\emptyset$, guaranteeing that our mapping preserves $E$. Moreover, for any fixed $W$ and two distinct $M,M'\in A_{i,X}\cap E$ their respective $\mathring{\Gamma}$ are distinct. Indeed, consider $\{\{v_\alpha,v_\beta\},\{v_\chi,v_\delta\}\}\subset M\setminus X$ and $\{\{v_\alpha,v_\chi\},\{v_\beta,v_{\delta'}\}\}\subset M'\setminus X$, where $\{v_{\alpha}, v_{\beta}\}$ and $\{v_{\alpha}, v_{\chi}\}$ is the first edge in which $M$, $M'$ differ (after performing $\tau$ on $M\setminus X$ and $M'\setminus X$). Without loss of generality, let $\alpha = \min\{\alpha, \beta, \chi, \delta, \delta'\}$. Suppose that $\alpha$ is within the first $2p$ vertices of the ordered $(n-x)$-tuple (otherwise we are done since $\{\{v_\alpha,v_\beta\},\{v_\chi,v_\delta\}\}\subset M_W$ and $\{\{v_\alpha,v_\chi\},\{v_\beta,v_\delta\}\}\subset M'_W$). Then, if the mapping yields, for example, $\{\{v_4,v_\alpha\},\{v_5,v_\beta\}\}\subset M_W$, then we would yield $\{\{v_4,v_\alpha\},\{v_5,v_\chi\}\}\subset M'_W$; clearly $M_W\neq M'_W$. Thus there is an injection from $A_{i,X}\cap E$ to $B_W\cap E$, i.e. $|B_W\cap E|\geq |A_{i,X}\cap E|$.\bigskip

Therefore, we have $$p_0=\mathbb P\left(A_{i,X}\mid E\right)\leq\mathbb P\left(B_W\mid E\right).$$\bigskip

Summing over all possible $W$, we have $p_0\cdot\sum_{k=2x-n}^{x}{{2x\choose 2k}\frac{(2k)!}{k!\cdot2^k}}\leq 1$, which gives $$p_0\leq\frac{1}{\sum_{k=2x-n}^{x}{{2x\choose 2k}\frac{(2k)!}{k!\cdot2^k}}}=p.$$\bigskip

Now, we want $ep(d+1)\leq 1$. This is equivalent to $$k\leq\frac{1}{e\cdot2x(2n-1){n-1\choose x-1}}\left(\sum_{j=2x-n}^{x}{{2x\choose 2j}\frac{(2j)!}{j!\cdot2^j}}-e\right).$$\smallskip

\begin{flushright}
$\blacksquare$
\end{flushright}
\end{shaded}


\section{A Conjecture}

In the proof of Theorem \ref{Main}, we used the idea of $\tau$ transformation to create unique permutations of perfect matchings. Here, we extend the notion of perfect matchings of graphs to that of $t$-uniform hypergraphs of order $tn$, i.e. $\forall v\in V\exists! e\in E(v\in e)$. This gives us a more general problem as follows:\bigskip

\begin{framed}
Let $\mathcal M=\{M_1,...,M_s\}$ be a collection of perfect matchings of a $t$-uniform hypergraph $H$ of order $tn$. Moreover, each of the $\frac{(nt)!}{(t!)^n\cdot n!}$ $t$-edges of the hypergraph appears at most $k$ times. Suppose that there does not exist a perfect matching of $H$ which agrees with each perfect matching $M_i\in\mathcal M$ in at most $x-1$ edges. What is the best possible lower bound for $k$?
\end{framed}

\bigskip

Following the method of proof of Theorem \ref{Main}, randomly pick a perfect matching $M$ from the set of all perfect matchings of $H$. Consider any $M_i\in \mathcal M$ and any $x$-matching $X\subset M_i$. Let $A_{i,X}$ be the event that $X\subset M$, and connect $A_{i,X}$ to $A_{i',X'}$ iff $X$ and $X'$ share at least one common vertex in their underlying vertex set. For each $A_{i,X}$, let the set of its neighbours in the dependency graph be $\mathcal D_{i,X}$. A combinatorial argument yields $$|\mathcal D_{i,X}|\leq k\cdot tx{tn-1\choose t-1}{n-1\choose x-1}=d.$$\bigskip

If we attempt to bound $$\mathbb P\left(A_{i,X}\mid\bigcap_{A_{i',X'}\in\mathcal S}{\overline{A_{i',X'}}}\right)=p_0$$

for any subset $\mathcal S\subset\mathcal A\setminus\mathcal D_{i,X}$, a difficulty arises if we mirror the mapping technique. We can still consider events similar to $W$ which split the underlying vertex set of $X$ into sets of $t$-edges, $(t-1)$-edges etc. and order them. Moreover, if a transformation similar to $\tau$ is performed on any matching $M\in A_{i,X}\cap E$, injectivity is still preserved. However, while it seems intuitively true that our $B_W\cap E$ should contain more elements that $A_{i,X}\cap E$, it is not as straightforward to map vertices into the respective $(t-\theta)$-edges, $t-1\geq\theta\geq 1$, such that injectivity is preserved. Hence, the problem remains open. Particularly, we conjecture:\bigskip

\begin{con}
Let $\mathcal M=\{M_1,...,M_s\}$ be a collection of perfect matchings of a $t$-uniform hypergraph $H$ of order $tn$. Moreover, each of the $\frac{(nt)!}{(t!)^n\cdot n!}$ $t$-edges of the hypergraph appears at most $k$ times. If $k\leq \frac{1}{e\cdot tx{tn-1\choose t-1}{n-1\choose x-1}}\left(N-e\right)$, where $$N=\sum{\left[\prod_{i=1}^{n}{{tx\choose ta_i}\frac{(ta_i)!}{(t!)^{a_i}(a_i)!}}\right]}$$\smallskip
is the sum over all vectors $(a_1,...,a_t)\in\left(\mathbb N\cup\{0\}\right)^t$ satisfying $ta_t+(t-1)a_{t-1}+...+a_1=tx$ and $a_t+a_{t-1}+...+a_1\leq n$, then there exists a perfect matching which agrees with each perfect matching $M_i\in\mathcal M$ in at most $x-1$ edges.
\end{con}

\bigskip

The upper bound for $k$ is based on the assumption that the intuition is correct.

\section{Conclusion}

It is unknown whether the upper bound obtained for $k$, namely\bigskip

$$\frac{1}{e\cdot2x(2n-1){n-1\choose x-1}}\left(\sum_{j=2x-n}^{x}{{2x\choose 2j}\frac{(2j)!}{j!\cdot2^j}}-e\right),$$

\medskip

is optimal, inasmuch as there is hitherto no research done in this area. However, it is possibly a fairly strong bound because the Lov\'{a}sz sieve is known to establish good bounds in problems.\bigskip
 
A possible continuation of our research is as follows: In \cite{ku}, a semi-random construction of a permutation code was given. In particular, using an analogue of Theorem \ref{Lovasz} for two events, an algorithm was formulated to construct a set of permutations in $S_n$ that is $<s$-intersecting in polynomial expected time. It would be possible to consider an analogue of the semi-random construction for collections of perfect matchings in $K_{2n}$.

\newpage


\end{document}